# The Complexity of Three-Way Statistical Tables

Jesus De Loera [*]    Shmuel Onn [†]

*Dedicated to Bernd Sturmfels on the occasion of his 40th birthday.*

**Abstract**

Multi-way tables with specified marginals arise in a variety of applications in statistics and operations research. We provide a comprehensive complexity classification of three fundamental computational problems on tables: existence, counting and entry-security.

One major outcome of our work is that each of the following problems is intractable already for "slim" 3-tables, with constant and smallest possible number 3 of rows: (1) deciding existence of 3-tables with given consistent 2-marginals; (2) counting all 3-tables with given 2-marginals; (3) finding whether an integer value is attained in entry $(i, j, k)$ by at least one of the 3-tables satisfying given (feasible) 2-marginals. This implies that a characterization of feasible marginals for such slim tables, sought by much recent research, is unlikely to exist.

Another important consequence of our study is a systematic efficient way of embedding the set of 3-tables satisfying any given 1-marginals and entry upper bounds in a set of slim 3-tables satisfying suitable 2-marginals with no entry bounds. This provides a valuable tool for studying multi-index transportation problems and multi-index transportation polytopes.

## 1 Introduction

A *d-table of size* $(n_1, \ldots, n_d)$ is an array of nonnegative integers $v = (v_{i_1,\ldots,i_d})$, $1 \leq i_j \leq n_j$. For $0 \leq m < d$, an *m-marginal* of $v$ is any of the $\binom{d}{m}$ possible $m$-tables obtained by summing the entries over all but $m$ indices. For instance, if $(v_{i,j,k})$ is a 3-table then its 0-marginal is $v_{+,+,+} = \sum_{i=1}^{n_1} \sum_{j=1}^{n_2} \sum_{k=1}^{n_3} v_{i,j,k}$, its 1-marginals are $(v_{i,+,+}) = (\sum_{j=1}^{n_2} \sum_{k=1}^{n_3} v_{i,j,k})$ and likewise $(v_{+,j,+})$, $(v_{+,+,k})$, and its 2-marginals are $(v_{i,j,+}) = (\sum_{k=1}^{n_3} v_{i,j,k})$ and likewise $(v_{i,+,k})$, $(v_{+,j,k})$.

Such tables appear naturally in statistics and operations research under various names such as *multi-way contingency tables*, *transportation matrices*, or *tabular data*. In all these applications, the tables of interest are those satisfying various constraints such as specified marginals or specified upper and lower bounds on the various table entries. Tables are central products of statistical agencies (for example, see the site [13] of the U.S Bureau of Census).

In this article we study three essential computational problems of constrained tables, primarily motivated by applications in statistical analysis and statistical data security (see e.g.

---

[*]Research partially supported by NSF grant 0073815
[†]Research partially supported by a grant from ISF - the Israel Science Foundation.

[9, 11, 14] and references therein): the **table existence or feasibility problem**, the **table counting problem** and the **table entry-security problem**. We provide comprehensive computational complexity classifications of these problems, which are discussed, respectively in Subsections 1.1-1.3 below, where we make the precise definition, include a briefing on the motivating statistical background, and describe our results for each problem. On the way we also show that the set of 3-tables satisfying any given 1-marginals and upper bounds on entries can be embedded in a set of "slim" 3-tables satisfying suitable 2-marginals with no entry bounds; this is discussed further, with the implications for the class of so-called *multi-index transportation polytopes* (cf. [23, 28]) in Subsection 1.4 below.

We demonstrate that even *slim* 3-tables, of size $(3, c, h)$, having fixed number 3 of rows, can have an arbitrarily complex behavior. This improves on earlier results of Irving and Jerrum [20]. For statisticians and agencies manipulating tabular data our results have practical repercussions: polynomial time algorithms for solving any one of the problems (feasibility, counting, or entry-security) are unlikely to exist. Thus, in taming 2-marginals arising in practice it will be necessary to exploit particular features of the real data in each specific application.

Our results on the intractability of slim 3-tables stand in contrast with the efficient methods available for 2-way contingency tables (see [16, 17]) and for the so-called decomposable-graph-log-linear-models [10, 14, 18]. Thus, we settle a problem that has been the focus of much recent research (see [6, 11, 25] and further references therein), on trying to find efficient methods for slim 3-tables, and demonstrate that already the 3-tables of smallest possible size $(3, c, h)$ which are not decomposable-graph-log-linear-models can have an arbitrarily complex behavior.

Finally, we point out that our results on the intractability of 2-marginals in 3-tables obviously extend to higher dimensions as 2-marginals in 3-tables can be embedded in higher dimensions.

## 1.1 Table existence or feasibility

First, we consider the **table existence problem**, also called the **feasibility problem** (cf. [4, 28]): *Given a prescribed collection of marginals that seem to describe a d-table of size $(n_1, \ldots, n_d)$, does there really exist a table with these marginals (and can it be effectively determined) ?* This problem is relevant for statistical analysis; for instance, disclosed or transmitted marginals may become perturbed or distorted in such a way that a feasible table may no longer exist, in which case, not only the data looses utility to the users, but also algorithms such as the iterative proportional fitting can fail to converge. This can be a problem because several statistical procedures are insensitive to existence, e.g. Fréchet type bounds presented in [14]. See [4, 5] for a discussion of the importance of table existence in statistics applications. An obvious necessary condition for the existence of a table with a specified collection of marginals is that the collection is *consistent*, that is, any two given marginals must agree on any of their common lower dimensional marginals. For instance, for the existence of a 3-table with specified 2-marginals $(v_{i,j,+})$ and $(v_{i,+,k})$, these marginals must agree on their common 1-marginal $(v_{i,+,+})$, so the 1-table equation $(\sum_{j=1}^n v_{i,j,+}) = (\sum_{k=1}^n v_{i,+,k})$ must hold. In general, however, these consistency equations do not even guarantee the existence of an array with nonnegative *real* entries (cf. [23, 28]).

The existence problem is easy to solve for 2-tables or 1-marginals, so the first really interesting case is that of 3-tables with all 2-marginals specified. The following theorem provides an almost complete classification of the complexity of this problem. We assume without loss of generality that the size $(r,c,h)$ of the tables satisfies $3 \leq r \leq c \leq h$: for $r \leq 2$ the problem reduces to the well-studied case of 2-dimensional tables and can be solved in polynomial time using linear programming over the corresponding *transportation polytope*: see further discussion of these polytopes and their higher dimensional *multi-index* generalizations in Subsection 1.4 below.

**Theorem 1.1** *The computational complexity of the existence problem for 3-tables of size $(r,c,h)$ with $3 \leq r \leq c \leq h$ and all 2-marginals specified is provided by the following table:*

|  | $r,c,h$ fixed | $r,c$ fixed, $h$ variable | $r$ fixed, $c,h$ variable | $r,c,h$ variable |
|---|---|---|---|---|
| *unary* 2-marginals | P | P | NPC | NPC |
| *binary* 2-marginals | P | ? | NPC | NPC |

Each entry $(i,j)$ of this table ($i=1,2$, $j=1,2,3,4$) represents a refined version of the existence problem; any problem $(2,j)$ is at least as hard as the problem $(1,j)$ below it, and for $j \geq 2$, any problem $(i,j)$ is at least as hard as the problem $(i,j-1)$ to its left. Here $NPC$ stands for *NP-complete* hence presumably intractable and practically unsolvable for large inputs (cf. [15]), whereas $P$ stands for *polynomial time* hence efficiently solvable. Binary versus unary are the two standard ways of encoding numbers: for *binary* marginals, the size of each marginal is taken to be the number of digits in its binary (or decimal) expansion hence is proportional to its logarithm, whereas for *unary* marginals, the size of each marginal is the marginal itself.

Thus, some of the entries follow at once from others, but are included for complete classification. Entry $(1,1)$ easily follows from exhaustive search. However, already entry $(2,1)$ requires the sophisticated algorithm of Lenstra for integer programming in fixed dimension [24]: it would be interesting to device a special faster polynomial time algorithm for this entry. In Section 3 we shall prove entry $(1,2)$, that is, the polynomial time solvability of existence for unary marginals with $r,c$ fixed ("small") and $h$ variable ("large"). In Section 2 we shall establish entries $(1,3)$ and $(2,3)$, that is, the NP-completeness of existence for marginals with $r=3$ fixed and $c,h$ variable. This implies at once entries $(1,4)$ and $(2,4)$ in the right column established previously by Irving and Jerrum [20], strengthening their results. Entry $(2,2)$ remains unsettled and challenging.

## 1.2 Table counting

Next, we consider the **table counting problem**: *given a prescribed collection of marginals, how many d-tables are there that share these marginals ?* Table counting has several applications in statistical analysis, in particular independence testing, and has been the focus of much research (see [8, 9, 21] and the extensive list of references therein). The counting problem can be formulated as that of counting the number of integer points in the associated multi-index transportation polytope (see further discussion in Subsection 1.4 below). The following analogue of Theorem 1.1 provides a complete classification of the complexity of this problem.

**Theorem 1.2** *The computational complexity of the counting problem for 3-tables of size $(r,c,h)$ with $2 \leq r \leq c \leq h$ and all 2-marginals specified is provided by the following table:*

|                    | $r,c,h$ fixed | $r,c$ fixed, $h$ variable | $r$ fixed, $c,h$ variable | $r,c,h$ variable |
|--------------------|---------------|---------------------------|---------------------------|------------------|
| *unary* 2-marginals  | $P$           | $P$                       | $\#PC$                    | $\#PC$           |
| *binary* 2-marginals | $P$           | $\#PC$                    | $\#PC$                    | $\#PC$           |

Here $\#PC$ stands for *$\#P$-complete*, hence presumably intractable; see Valiant's seminal paper [27] which introduces the complexity theory of counting, or consult [15].

Entry $(1,1)$ is easy, but already entry $(2,1)$ requires the sophisticated algorithm of Barvinok for counting integer points in polytopes in fixed dimension [3]. In Section 3 we prove entry $(1,2)$, that is, counting in polynomial time for unary marginals with $r,c$ fixed and $h$ variable. Entry $(2,2)$ follows from the $\#P$-completeness of counting 2-tables of size $(2,n)$ with binary marginals [12]. In Section 3 we also prove entry $(1,3)$, that is, the $\#P$-completeness of counting for unary marginals with $r=2$ fixed and $c,h$ variable, implying entries $(2,3)$, $(1,4)$ and $(2,4)$ as well.

### 1.3 Table entry-security

The third problem we consider arises in the context of secure and confidential disclosure of public statistical data (see [6, 11, 14] and references therein). The goal in this context is the release of some marginals of a table in the database but not the table's entries themselves. If the range of possible values that an entry $x_{i,j,k}$ can attain in any table satisfying the released collection of marginals is too narrow, or even worse, consists of the unique value of that entry in the actual table in the database, then this entry may be exposed. This shows the importance of determining tight integer upper and lower bounds of each entry $x_{i,j,k}$. We consider the following *lower* and *upper* versions of this problem, the **table entry-security problem**: *let there be given a* **feasible** *prescribed collection of 2-marginals (i.e., admitting at least one feasible d-table) and an index tuple $(i_1,\ldots,i_d)$. Given now also a nonnegative integer $L$, is there a d-table $x$ having the given marginals whose entry $x_{i_1,\ldots,i_d}$ is greater than or equal to $L$ (lower version) ? or the analogous question, given now also a nonnegative integer $U$, is there a d-table $x$ having the given marginals whose entry $x_{i_1,\ldots,i_d}$ is less than or equal to $U$ (upper version) ?*

There is an extensive work on the entry-security problem, see e.g. [5, 10, 14, 25, 26], where properties are sought that may help address the problem. The paper [11] surveys the state-of-the-art research and practical techniques. The following analogue of Theorems 1.1 and 1.2 provides an almost complete classification of the complexity of the entry-security problem as well.

**Theorem 1.3** *The complexity of both lower and upper versions of the entry-security problem for 3-tables of size $(r,c,h)$ with $4 \leq r \leq c \leq h$ and all 2-marginals specified is provided by:*

|                    | $r,c,h$ fixed | $r,c$ fixed, $h$ variable | $r$ fixed, $c,h$ variable | $r,c,h$ variable |
|--------------------|---------------|---------------------------|---------------------------|------------------|
| *unary* 2-marginals  | $P$           | $P$                       | $NPC$                     | $NPC$            |
| *binary* 2-marginals | $P$           | ?                         | $NPC$                     | $NPC$            |

Once again, entry $(1,1)$ is easy and entry $(2,1)$ follows from [24]. In Section 4 we prove the polynomial time solvability for the case of unary marginals with $r,c$ fixed and $h$ variable (entry $(1,2)$, and the intractability for $r$ fixed and $c,h$ variable (entries $(1,3)$ and $(2,3)$, strengthening earlier hardness results of Irving and Jerrum [20] reflected in entries $(1,4)$ and $(2,4)$. Again, entry $(2,2)$ remains unsettled.

## 1.4 Multi-index transportation polytopes and the power of two-marginals

Given a specified collection of marginals for $d$-tables of size $(n_1, \ldots, n_d)$, possibly together with specified lower and upper bounds on some of the table entries, the associated *multi-index transportation polytope* is the set of all nonnegative *real valued* arrays satisfying these marginals and entry bounds (cf. [23]), and is a (typically bounded) convex polyhedron in $\mathbf{R}^{n_1 \cdots n_d}$. For instance, for 2-tables of size $(n, n)$ with all 1-marginals equal to 1 and no entry bounds, this is the Birkhoff polytope of $n$ by $n$ bistochastic matrices. The $d$-tables satisfying the given marginals and entry bounds are precisely the *integer points* in the associated multi-index transportation polytope.

In Section 5 we show how a system of 1-marginal and entry upper bound constraints on 3-tables can be embedded into a system of 2-marginal constraints (with no entry bounds) on "slim" 3-tables, demonstrating the expressive power of 2-marginals and reducing the existence, counting and entry-security problems for 1-marginals with upper bounds to that for 2-marginals with no upper bounds in slim tables. We prove the following somewhat technical statement.

**Theorem 1.4** *Given 1-marginals $(u_{i,+,+})$, $(u_{+,j,+})$, $(u_{+,+,k})$ and entry upper bounds $(p_{i,j,k})$ for 3-tables of size $(r, c, h)$, there exist polynomial time constructible 2-marginals $(v_{i,j,+})$, $(v_{i,+,k})$, $(v_{+,j,k})$ for 3-tables of size $(3, rc, r + c + h)$ such that the set of nonnegative real $(r, c, h)$-arrays with the given upper bounds and 1-marginals is in integer preserving affine bijection with the set of nonnegative real $(3, rc, r + c + h)$-arrays with the constructed 2-marginals.*

A particularly appealing outcome of our constructions is the systematic derivation of "real-feasible-integer-infeasible" collections of 2-marginals, admitting nonnegative real 3-arrays but no (integer) 3-tables. Remarkably, we can "automatically" obtain the 2-marginals for 3-tables of size $(3, 4, 6)$, discovered by Vlach [28, page 77] (see our Example 2.2). We obtain it from very simple $\{0, 1\}$-valued 1-marginals and entry upper bounds for 3-tables of size $(2, 2, 2)$. Applying to them the constructions of this paper yields precisely the collection of $\{0, 1\}$-valued 2-marginals discovered by Vlach.

We conclude this introduction with some final discussion. First, we refer to the open problems left in the $(2, 2)$ entries of Theorems 1.1 and 1.3. Consider the set of 3-tables $v$ of size $(r, c, h)$ with $r, c$ fixed satisfying specified marginals $(v_{i,+,k})$ and $(v_{+,j,k})$ but *without* restriction on the marginal $(v_{i,j,+})$. The projection $\mathbf{R}^{r \cdot c \cdot h} \longrightarrow \mathbf{R}^{r \cdot c} : v \mapsto (v_{i,j,+}) = \sum_{k=1}^{h} v_{i,j,k}$ sends the associated multi-index transportation polytope onto a subpolytope $P$ of the transportation polytope of all 2-tables of size $(r, c)$ with 1-marginals $(u_{i,+}) := (v_{i,+,+})$ and $(u_{+,j}) := (v_{+,j,+})$. The techniques of [1, 2, 19, 22] allow to produce the vertices of $P$ in polynomial time and check if any given "vertical" marginal $(v_{i,j,+})$ lies in $P$, which is a necessary condition for the existence of a 3-table with $(v_{i,+,k})$, $(v_{+,j,k})$ and $(v_{i,j,+})$. Further development of the methods of [1, 2, 19, 22] combined with integer programming in fixed dimension might help in addressing these remaining problems.

Although our results stress the complexity of handling even small and slim 3-way tables for statistical applications, recent results using special structure in specific systems may make such table systems amenable to geometric algorithms for practical computations. For example, in [14, 10, 18] the specified marginals satisfies a hierarchical structure of certain graphical models in statistics. Other approaches include the new generation of algebraic and randomized algorithms [7, 12], which will allow, in practice, faster computations for increasingly larger problems.

## 2   The table existence problem

In this section we provide the proof of Theorem 1.1 discussed in Subsection 1.1, and demonstrate our constructions with some examples. In particular, as explained in the introduction, using our construction we recover the smallest possible real-feasible-integer-infeasible collection of 2-marginals of 3-tables of size $(3, 4, 6)$.

**Proof of Theorem 1.1.** As explained in Subsection 1.1, entry $(1,1)$ of the complexity table claimed by the Theorem is easy and entry $(2,1)$ follows from [24]. Entry $(1,2)$ follows from entry $(1,2)$ in the table of Theorem 1.2 which will be proved in the next section: indeed, we shall show in Section 3 how to compute in polynomial time the number of 3-tables of size $(r, c, h)$ with $r, c$ fixed satisfying given 2-marginals in unary, and hence in particular to decide if this number is zero or not, providing a solution of the existence problem as well.

We need then prove entry $(1,3)$ of the table, which implies at once entries $(2,3), (1,4), (2,4)$ as well. It is easy to see that the well-known 3-dimensional matching problem (cf. [15]) is equivalent to the following problem: given a $\{0,1\}$-valued 3-table $p = (p_{i,j,k})$ of size $(n, n, n)$, is there a 3-table $x = (x_{i,j,k})$ with all 1-marginals equal to 1 which is *dominated* by $p$, i.e. satisfies the upper bounds $x_{i,j,k} \leq p_{i,j,k}$ for all $i, j, k$ ? we reduce this problem to ours, which is clearly in NP. Let then $p = (p_{i,j,k})$ be a given $\{0,1\}$-valued 3-table of size $(n, n, n)$. We define efficiently constructible 2-marginals for $(3, n^2, 3n)$-tables such that nonnegative real arrays $y$ with these marginals are in integer preserving affine bijection with nonnegative real $(n, n, n)$-arrays $x$ with all 1-marginals equal to 1 dominated by $p$. For clarity, the table will be indexed by triplets of special form which we now explain. The first index will be an integer $1 \leq t \leq 3$. The second index will be an ordered pair $ij$ with $1 \leq i, j \leq n$. The third index will belong into one of three groups - "domination" group, "row" group and "column" group, and will consist of one of three three-letter abbreviations gro $\in \{\text{dom}, \text{row}, \text{col}\}$ according to the group it belongs to, along with a numerical index $1 \leq k \leq n$. The 2-marginals are provided by the following three matrices.

$$
(v_{+,ij,\text{gro }k}) = \begin{pmatrix}
\begin{array}{c|cccccccccccc}
 & 11 & 12 & \cdots & 1n & 21 & 22 & \cdots & 2n & \cdots & n1 & n2 & \cdots & nn \\
\hline
\text{dom } 1 & p_{1,1,1} & p_{1,2,1} & \cdots & p_{1,n,1} & p_{2,1,1} & p_{2,2,1} & \cdots & p_{2,n,1} & \cdots & p_{n,1,1} & p_{n,2,1} & \cdots & p_{n,n,1} \\
\text{dom } 2 & p_{1,1,2} & p_{1,2,2} & \cdots & p_{1,n,2} & p_{2,1,2} & p_{2,2,2} & \cdots & p_{2,n,2} & \cdots & p_{n,1,2} & p_{n,2,2} & \cdots & p_{n,n,2} \\
\cdots & \vdots & \vdots & \vdots & \vdots & \vdots & \vdots & \vdots & \vdots & \vdots & \vdots & \vdots & \vdots & \vdots \\
\text{dom } n & p_{1,1,n} & p_{1,2,n} & \cdots & p_{1,n,n} & p_{2,1,n} & p_{2,2,n} & \cdots & p_{2,n,n} & \cdots & p_{n,1,n} & p_{n,2,n} & \cdots & p_{n,n,n} \\
\hline
\text{row } 1 & 1 & 1 & \cdots & 1 & 0 & 0 & \cdots & 0 & \cdots & 0 & 0 & \cdots & 0 \\
\text{row } 2 & 0 & 0 & \cdots & 0 & 1 & 1 & \cdots & 1 & \cdots & 0 & 0 & \cdots & 0 \\
\cdots & \vdots & \vdots & \vdots & \vdots & \vdots & \vdots & \vdots & \vdots & \vdots & \vdots & \vdots & \vdots & \vdots \\
\text{row } n & 0 & 0 & \cdots & 0 & 0 & 0 & \cdots & 0 & \cdots & 1 & 1 & \cdots & 1 \\
\hline
\text{col } 1 & 1 & 0 & \cdots & 0 & 1 & 0 & \cdots & 0 & \cdots & 1 & 0 & \cdots & 0 \\
\text{col } 2 & 0 & 1 & \cdots & 0 & 0 & 1 & \cdots & 0 & \cdots & 0 & 1 & \cdots & 0 \\
\cdots & \vdots & \vdots & \vdots & \vdots & \vdots & \vdots & \vdots & \vdots & \vdots & \vdots & \vdots & \vdots & \vdots \\
\text{col } n & 0 & 0 & \cdots & 1 & 0 & 0 & \cdots & 1 & \cdots & 0 & 0 & \cdots & 1 \\
\end{array}
\end{pmatrix}
$$

$$(v_{t,ij,+}) \;=\; \begin{array}{c} \\ 1 \\ 2 \\ 3 \end{array}\!\!\left(\begin{array}{ccccccccccccc} 11 & 12 & \cdots & 1n & 21 & 22 & \cdots & 2n & \cdots & n1 & n2 & \cdots & nn \\ \hline 1 & 1 & \cdots & 1 & 1 & 1 & \cdots & 1 & \cdots & 1 & 1 & \cdots & 1 \\ p_{1,1,+} & p_{1,2,+} & \cdots & p_{1,n,+} & p_{2,1,+} & p_{2,2,+} & \cdots & p_{2,n,+} & \cdots & p_{n,1,+} & p_{n,2,+} & \cdots & p_{n,n,+} \\ 1 & 1 & \cdots & 1 & 1 & 1 & \cdots & 1 & \cdots & 1 & 1 & \cdots & 1 \end{array}\right)$$

$$(v_{t,+,\text{gro}\,k}) \;=\; \begin{array}{c} \\ \text{dom } 1 \\ \text{dom } 2 \\ \cdots \\ \text{dom } n \\ \\ \text{row } 1 \\ \text{row } 2 \\ \cdots \\ \text{row } n \\ \\ \text{col } 1 \\ \text{col } 2 \\ \cdots \\ \text{col } n \end{array}\!\!\left(\begin{array}{ccc} 1 & 2 & 3 \\ \hline 1 & p_{+,+,1}-1 & 0 \\ 1 & p_{+,+,2}-1 & 0 \\ \vdots & \vdots & \vdots \\ 1 & p_{+,+,n}-1 & 0 \\ \hline n-1 & 0 & 1 \\ n-1 & 0 & 1 \\ \vdots & \vdots & \vdots \\ n-1 & 0 & 1 \\ \hline 0 & 1 & n-1 \\ 0 & 1 & n-1 \\ \vdots & \vdots & \vdots \\ 0 & 1 & n-1 \end{array}\right)$$

First consider any nonnegative real $(n, n, n)$-array $x$ dominated by $p$ with all 1-marginals equal 1. We show it uniquely extends to a nonnegative real $(3, n^2, 3n)$-array with 2-marginals as above. It will be convenient to keep in mind a partition of $(3, n^2, 3n)$-arrays into blocks as in Figure 1.

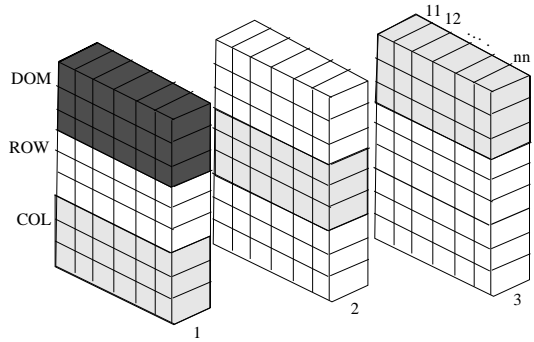

Figure 1: Block partition of $(3, n^2, 3n)$-arrays.

Given then such an array $x$, embed it in the **black block** $(1, \text{dom})$ of a $(3, n^2, 3n)$-array $y$ by

$$y_{1,ij,\text{dom}\,k} \;:=\; x_{i,j,k}, \qquad 1 \leq i,j,k \leq n \;.$$

We now show that the block $x$ can be uniquely extended to a whole nonnegative real $(3, n^2, 3n)$-array $y$ with the above 2-marginals. First consider the entries in the **grey blocks** $(1, \text{col})$, $(2, \text{row})$ and $(3, \text{dom})$ in Figure 1: since $v_{1,+,\text{col } k} = v_{2,+,\text{row } k} = v_{3,+,\text{dom } k} = 0$ for all $k$, it follows that all the entries $y_{1,ij,\text{col } k}$, $y_{2,ij,\text{row } k}$ and $y_{3,ij,\text{dom } k}$ constituting these blocks are zero. Next, consider the entries in the **white block** $(1, \text{row})$: using the fact just established that all entries in the block $(1, \text{col})$ below it are zero, and examining the 2-marginals $v_{+,ij,\text{row } k}$ and $v_{1,ij,+} = 1$, we find that $y_{1,ij,\text{row } i} = 1 - \sum_{k=1}^{n} y_{1,ij,\text{dom } k} = 1 - x_{i,j,+} \geq 0$ whereas for $k \neq i$ we have $y_{1,ij,\text{row } k} = 0$. This also yields the entries in the **white block** $(3, \text{row})$: we have $y_{3,ij,\text{row } i} = 1 - y_{1,ij,\text{row } i} = x_{i,j,+} \geq 0$ whereas for $k \neq i$ we have $y_{3,ij,\text{row } k} = 0$. Next, consider the entries in the **white block** $(2, \text{dom})$: using the fact that all entries in the block $(3, \text{dom})$ to its right are zero, and examining the 2-marginals $v_{+,ij,\text{dom } k} = p_{i,j,k}$ we find that $y_{2,ij,\text{dom } k} = p_{i,j,k} - x_{i,j,k} \geq 0$ for all $i,j,k$. Next consider the entries in the **white block** $(2, \text{col})$: using the fact that all entries in the block $(2, \text{row})$ above it are zero, and examining the 2-marginals $v_{+,ij,\text{col } k}$ and $v_{2,ij,+} = p_{i,j,+}$, we find that $y_{2,ij,\text{col } j} = p_{i,j,+} - \sum_{k=1}^{n} y_{2,ij,\text{dom } k} = x_{i,j,+} \geq 0$ whereas for $k \neq j$ we have $y_{2,ij,\text{col } k} = 0$. This also yields the entries in the **white block** $(3, \text{col})$: we have $y_{3,ij,\text{col } j} = 1 - y_{2,ij,\text{col } j} = 1 - x_{i,j,+} \geq 0$ whereas for $k \neq j$ we have $y_{3,ij,\text{col } k} = 0$.

Next consider any nonnegative real $(3, n^2, 3n)$-array $y$ with the above 2-marginals, and let $x$ be its $(n,n,n)$-subarray given by the black block $(1, \text{dom})$ of $y$, defined by $x_{i,j,k} := y_{1,ij,\text{dom } k}$ for all $i,j,k$. We show that $x$ is nonnegative, dominated by $p$ and has all 1-marginals equal to 1. It is nonnegative since so is $y$. It is dominated by $p$ since, for all $i,j,k$ we have $p_{i,j,k} - x_{i,j,k} = y_{2,ij,\text{dom } k} \geq 0$. Finally, all the 1-marginals of $x$ are equal to 1 since:

$$x_{+,+,k} = \sum_{i,j} y_{1,ij,\text{dom } k} = v_{1,+,\text{dom } k} = 1, \qquad 1 \leq k \leq n\ ;$$

$$x_{i,+,+} = \sum_{j} y_{3,ij,\text{row } i} = v_{3,+,\text{row } i} = 1, \qquad 1 \leq i \leq n\ ;$$

$$x_{+,j,+} = \sum_{i} y_{2,ij,\text{col } j} = v_{2,+,\text{col } j} = 1, \qquad 1 \leq j \leq n\ .$$

Thus, the set of nonnegative real $(n,n,n)$-arrays $x$ dominated by $p$ and with all 1-marginals 1 is in integer preserving affine bijection with the set of nonnegative real $(3, n^2, 3n)$-arrays $y$ with the constructed 2-marginals. In particular, the corresponding sets of tables are in bijection and therefore the former is nonempty if and only if the latter is. This completes the reduction of 3-dimensional matching to our problem, the proof of entry $(1,3)$ and the proof of Theorem 1.1. □

The following two examples illustrate our construction.

**Example 2.1** Let $n = 2$ and let $p$ be the $\{0,1\}$ valued 3-table of size $(2,2,2)$ given by:

$$p_{1,1,1} = 1,\ p_{1,2,1} = 1,\ p_{2,1,1} = 0,\ p_{2,2,1} = 0,\ p_{1,1,2} = 1,\ p_{1,2,2} = 1,\ p_{2,1,2} = 0,\ p_{2,2,2} = 1\ .$$

Our construction yields the 2-marginals for 3-tables of size $(3,4,6)$ presented in Figure 2.

The unique 3-table $x$ with all 1-marginals equal to 1 which is dominated by $p$ is given by:

$$x_{1,1,1} = 1,\ x_{1,2,1} = 0,\ x_{2,1,1} = 0,\ x_{2,2,1} = 0,\ x_{1,1,2} = 0,\ x_{1,2,2} = 0,\ x_{2,1,2} = 0,\ x_{2,2,2} = 1,$$

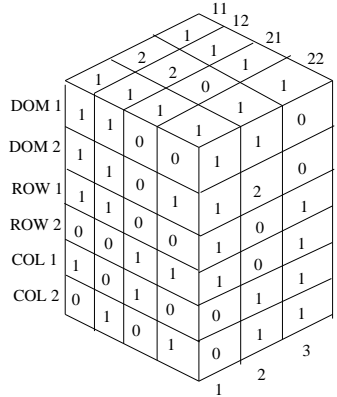

Figure 2: The 2-marginals $v_{+,ij,\text{gro } k}$, $v_{t,+,\text{gro } k}$, $v_{t,ij,+}$ constructed in Example 2.1.

and the corresponding 3-table $y$ with the above 2-marginals is given by the following blocks:

$$
\begin{array}{c}
y_{1,ij,\text{gro } k} \\
\begin{array}{c|cccc}
 & 11 & 12 & 21 & 22 \\
\hline
\text{dom } 1 & 1 & 0 & 0 & 0 \\
\text{dom } 2 & 0 & 0 & 0 & 1 \\
\hline
\text{row } 1 & 0 & 1 & 0 & 0 \\
\text{row } 2 & 0 & 0 & 1 & 0 \\
\hline
\text{col } 1 & 0 & 0 & 0 & 0 \\
\text{col } 2 & 0 & 0 & 0 & 0
\end{array}
\end{array}
,\quad
\begin{array}{c}
y_{2,ij,\text{gro } k} \\
\begin{array}{c|cccc}
 & 11 & 12 & 21 & 22 \\
\hline
\text{dom } 1 & 0 & 1 & 0 & 0 \\
\text{dom } 2 & 1 & 1 & 0 & 0 \\
\hline
\text{row } 1 & 0 & 0 & 0 & 0 \\
\text{row } 2 & 0 & 0 & 0 & 0 \\
\hline
\text{col } 1 & 1 & 0 & 0 & 0 \\
\text{col } 2 & 0 & 0 & 0 & 1
\end{array}
\end{array}
,\quad
\begin{array}{c}
y_{3,ij,\text{gro } k} \\
\begin{array}{c|cccc}
 & 11 & 12 & 21 & 22 \\
\hline
\text{dom } 1 & 0 & 0 & 0 & 0 \\
\text{dom } 2 & 0 & 0 & 0 & 0 \\
\hline
\text{row } 1 & 1 & 0 & 0 & 0 \\
\text{row } 2 & 0 & 0 & 0 & 1 \\
\hline
\text{col } 1 & 0 & 0 & 1 & 0 \\
\text{col } 2 & 0 & 1 & 0 & 0
\end{array}
\end{array}.
$$

As pointed out in the introduction, our construction can be used to systematically obtain "real-nonempty-integer-empty" multi-index transportation polytopes, namely, collections of 2-marginals admitting nonnegative real 3-arrays but no (integer) 3-tables. In particular, we next recover the smallest such example, first discovered by Vlach [28], as follows.

**Example 2.2** Let again $n = 2$ and let $p$ be the $\{0, 1\}$ valued 3-table of size $(2, 2, 2)$ given by:

$$p_{1,1,1} = 1, \ p_{1,2,1} = 0, \ p_{2,1,1} = 0, \ p_{2,2,1} = 1, \ p_{1,1,2} = 0, \ p_{1,2,2} = 1, \ p_{2,1,2} = 1, \ p_{2,2,2} = 0 \ .$$

Our construction yields the 2-marginals for 3-tables of size $(3, 4, 6)$ presented in Figure 3.

It can be verified that there is a single nonnegative real array of size $(2, 2, 2)$ with all 1-marginals equal to 1 which is dominated by the upper-bound table $p$. All entries of this array are $\{0, \frac{1}{2}\}$-valued and there is no (integer) table with the prescribed constraints. Our construction lifts this situation to 2-marginals with no upper bounds: all entries of the unique corresponding nonnegative real array of size $(3, 4, 6)$ with the 2-marginals in Figure 3 are $\{0, \frac{1}{2}\}$-valued and there is no (integer) table with these constructed 2-marginals.

Figure 3: Derivation of Vlach's example from our construction.

## 3   The table counting problem

In this section we provide the proof of Theorem 1.2 discussed in Subsection 1.2.

**Proof of Theorem 1.2.** As explained in Subsection 1.2, entry $(1,1)$ of the complexity table claimed by the theorem is easy, entry $(2,1)$ follows from [3], and entry $(2,2)$ follows from [12].

First, we prove entry $(1,3)$ of the table, which implies at once entries $(2,3),(1,4),(2,4)$ as well. We describe a direct reduction from Valiant's canonical $\#P$-complete problem of computing the *permanent* of a $\{0,1\}$-valued matrix [27] (recall that the *permanent* of an $n$ by $n$ matrix $A$ is $\mathrm{perm}(A) := \sum_\sigma \prod_{i=1}^n A_{i,\sigma(i)}$, the sum extending over all permutations $\sigma$ of $\{1,\ldots,n\}$; for instance, the permanent of the adjacency matrix of a subgraph of the complete bipartite graph $K_{n,n}$ is the number of perfect matchings in that subgraph).

Let then $A$ be a $\{0,1\}$-valued $n$ by $n$ matrix the permanent of which is to be computed. Define 2-marginals for 3-tables of size $(2,n,n)$ by

$$v_{i,j,+} := A_{i,j}, \quad v_{i,+,1} := v_{+,j,1} := 1, \quad v_{i,+,2} := A_{i,+} - 1, \quad v_{+,j,2} := A_{+,j} - 1, \qquad 1 \leq i,j \leq n.$$

Any 3-table $x$ with these marginals is determined by its 2-subtable $(x_{i,j,1})$ since for all $i,j$ we have $x_{i,j,2} = A_{i,j} - x_{i,j,1}$. Now, it is not hard to see that a nonnegative integer $n$ by $n$ matrix $\Sigma$ can arise as the subtable $(x_{i,j,1})$ of a 3-table $x$ with the constructed 2-marginals if and only if it is the standard representing matrix of a permutation $\sigma$ satisfying $\prod_{i=1}^n A_{i,\sigma(i)} = 1$. Therefore, the permanent of $A$, which is the number of such permutations $\sigma$, is precisely the number of 3-tables with the constructed 2-marginals, completing the reduction and the proof of entry $(1,3)$.

Next we prove entry $(1,2)$ of the table of Theorem 1.2. Note that, as explained in the proof of Theorem 1.1, this implies the corresponding entry $(1,2)$ in the table of Theorem 1.1 as well.

So, $r,c$ are fixed and we are given unary presented 2-marginals $v_{i,j,+}$, $v_{i,+,k}$ and $v_{+,j,k}$ for 3-tables of size $(r,c,h)$. Let $S$ be the set of all 2-tables $s$ of size $(r,c)$ satisfying the upper bounds $s_{i,j} \leq v_{i,j,+}$ for all $i,j$, that is, dominated by the given "vertical" marginals. For $k=1,\ldots,h$ define a matrix $A_k$ whose rows and columns are indexed by the elements of $S$, with entries

$$(A_k)_{s,t} := \begin{cases} 1 & \text{if } (t-s)_{i,+} = v_{i,+,k} \text{ for all } i \text{ and } (t-s)_{+,j} = v_{+,j,k} \text{ for all } j \\ 0 & \text{otherwise} \end{cases} \quad s,t \in S. \quad (1)$$

For $p = 1, \ldots, h$ let $A^p := A_1 \cdot A_2 \cdot \ldots \cdot A_p$ be the product of the matrices $A_k$, $k = 1, \ldots, p$. Let further $l, u$ denote, respectively, the tables in $S$ with entries $l_{i,j} := 0$ and $u_{i,j} := v_{i,j,+}$ for all $i, j$.

We claim that for any $1 \leq p \leq h$ and for any $s, t \in S$, the number of 3-tables $x$ of size $(r, c, p)$ with $x_{i,+,k} = v_{i,+,k}$, $x_{+,j,k} = v_{+,j,k}$ and $x_{i,j,+} = (t-s)_{i,j}$ for $1 \leq i \leq r$, $1 \leq j \leq c$ and $1 \leq k \leq p$, is precisely equal to the entry $A^p_{s,t}$ of $A^p$. In particular, the number of $(r, c, h)$-tables with the given 2-marginals is given by $A^h_{l,u}$. Since $r, c$ are fixed and the 2-marginals are presented in unary, the number $\prod_{i=1}^r \prod_{j=1}^c (v_{i,j,+} + 1)$ of tables in $S$ is polynomial in the size of the input and therefore the matrix $A^h$ and its sought entry $A^h_{l,u}$ can be computed in polynomial time.

We prove the claim by induction on $p$. First, consider the case $p = 1$ and let $s, t$ be any pair of tables of $S$. There is a unique $(r, c, 1)$-array $x$ satisfying $x_{i,j,1} = x_{i,j,+} = (t-s)_{i,j}$ for all $i, j$, and $x$ is a table satisfying $x_{i,+,1} = v_{i,+,1}$ and $x_{+,j,1} = v_{+,j,1}$ if and only if $(t-s)_{i,+} = v_{i,+,1}$ and $(t-s)_{+,j} = v_{+,j,1}$ for all $i, j$, which by Equation 1 holds if and only if $A^1_{s,t} = (A_1)_{s,t} = 1$.

Next, consider any $2 \leq p \leq h$ and suppose the statement is true for all values less than $p$. Let $s, t$ be any pair of tables of $S$. Then any $(r, c, p)$-table $x$ with $x_{i,+,k} = v_{i,+,k}$, $x_{+,j,k} = v_{+,j,k}$ and $x_{i,j,+} = (t-s)_{i,j}$ for all $i, j, k$ is obtained, for some $w \in S$, by augmenting any of the $A^{p-1}_{s,w}$ tables $y$ of size $(r, c, p-1)$ with $y_{i,+,k} = v_{i,+,k}$, $y_{+,j,k} = v_{+,j,k}$ and $y_{i,j,+} = (w-s)_{i,j}$ for all $i, j, k$ by any of the $(A_p)_{w,t}$ tables $z$ of size $(r, c, 1)$ with $z_{i,+,1} = v_{i,+,p}$, $z_{+,j,1} = v_{+,j,p}$ and $z_{i,j,+} = (t-w)_{i,j}$ for all $i, j$. Thus, the number of such tables is $\sum_{w \in S} A^{p-1}_{s,w} (A_p)_{w,t}$ which is precisely $A^p_{s,t}$, proving the induction step and the claim, thus completing the proof of entry $(1, 2)$ and of Theorem 1.2. □

## 4  The table entry-security problem

In this section we provide the proofs of Theorems 1.3 discussed in Subsections 1.3.

**Proof of Theorem 1.3.** Once again, as explained in Subsection 1.3, entry $(1, 1)$ of the complexity table claimed by the theorem is easy and entry $(2, 1)$ follows from [24]. We need to prove entry $(1, 3)$ which implies at once entries $(2, 3), (1, 4), (2, 4)$ as well, and entry $(1, 2)$.

We begin with the proof of entry $(1, 3)$: we prove that both the lower and upper versions of the entry-security problem are hard. In fact, we prove stronger results by showing that each of the following special cases are already hard: (A) we reduce the table feasibility problem to the problem of deciding whether, given a feasible collection of 2-marginals, there is a 3-table with a specified entry equal to the minimal possible value zero; thus, even the special case of the upper version of the entry security-problem with $U = 0$ is hard; (B) we reduce 3-dimensional matching to the problem of deciding whether, given a feasible collection of 2-marginals, there is a 3-table with a specified entry equal to the maximal possible value given by the Fréchet upper bound (minimal value of the three 2-marginals involving this entry). Thus, even the special case of the lower version of the entry-security-problem with $L$ the Fréchet upper bound is hard.

We being with part (A). Suppose then that we are given all the 2-marginals $v_{i,j,+}$, $v_{+,j,k}$, $v_{i,+,k}$ for 3-tables of size $(r, c, h)$, and we wish to know whether there is indeed a 3-table with these 2-marginals. Observe that the input 2-marginals can be assumed to be consistent; otherwise the 2-marginals are infeasible. The 2-marginals given are $A_{i,j} = v_{i,j,+}, B_{j,k} = v_{+,j,k}$ and $C_{i,k} = v_{i,+,k}$. It is important to observe that $\sum_i^r \sum_k^h C_{i,k} = \sum_i^r \sum_j^c A_{i,j} = \sum_j^c \sum_k^h B_{j,k} = T$. Note that $T$ denotes the total sum of all entries on any 3-table that satisfies the input 2-marginals $A, B, C$. They also share the same 1-marginals $C_{i,+} = A_{i,+}, A_{+,j} = B_{j,+}$, and $C_{+,k} = B_{j,+}$. All these equalities follow because the 2-marginals are consistent and will be useful later on.

Now we will construct a *feasible* set of 2-marginals for a family of 3-tables $R_{s,t,u}$ of size $(r+1, c+1, h+1)$. The entry-value of a certain entry $R_{s,t,u}$ can be used to decide whether the original set of 2-marginals $A, B, C$ is feasible. We present the 2-marginals in Figure 4 as numbers on the surface of a 3-table. The three 2-marginals $R_{+,t,u}$, $R_{s,+,u}$, $R_{s,t,+}$ are indicated by the coordinate directions in Figure 4. The reader can verify (see Figure 4) that the assignment is done as follows: for 2-marginal $R_{s,+,u}$ we set $R_{1,+,1} = T$, $R_{1,+,2} = C_{+,h}, R_{1,+,3} = C_{+,h-1}, \ldots R_{1,+,t} = C_{+,h-t+2}, \ldots, R_{1,+,h+1} = C_{+,1}$. Similarly $R_{2,+,1} = C_{1,+}, \ldots, R_{s,+,1} = C_{s-1,+}, \ldots R_{r+1,+,1} = C_{r,+}$. Finally, we have the assignment $R_{s,+,u} = C_{s-1,u-1}$ for $s = 2 \ldots r+1$ and $u = 2 \ldots h+1$. Next for 2-marginal $R_{s,t,+}$ we have $R_{1,1,+} = T, R_{1,2,+} = B_{1,+}, \ldots, R_{1,t,+} = B_{t-1,+}, \ldots R_{1,c+1,+} = B_{c,+}$, we also have $R_{2,1,+} = C_{1,+}, \ldots, R_{s,1,+} = C_{s-1,+}, \ldots R_{r+1,1,+} = C_{r,+}$, and $R_{s,t,+} = A_{s-1,t-1}$ for $s = 2 \ldots r+1$ and $t = 2, \ldots c+1$. Finally for 2-marginal $R_{+,t,u}$ we have that $R_{+,1,1} = T$ and we set $R_{+,2,1} = A_{+,1}, \ldots R_{+,t,1} = A_{+,t-1}, \ldots R_{+,c+1,1} = A_{+,c}$. Also from the picture we see that $R_{+,1,2} = C_{+,h}, \ldots R_{+,1,u} = C_{+,h-u+2} \ldots R_{+,1,h+1} = C_{+,1}$ and $R_{+,t,u} = B_{t-1,u-1}$ for $t = 2 \ldots c+1$ and $u = 2, \ldots, h+1$.

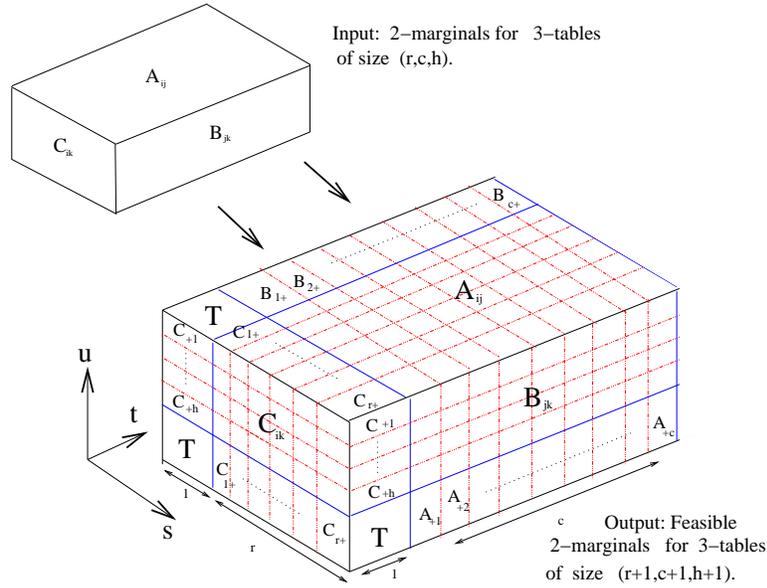

Figure 4: The Construction of feasible 2-marginals from input 2-marginals.

Note that any such 3-table $R$ with the 2-marginals we wrote breaks up naturally into eight smaller 3-tables. We show the blocks $b(1), \ldots, b(8)$ in Figure 5 marking their dimensions. Observe that block $b(i)$ has a set of 2-marginals corresponding to it, i.e. those that come from "projecting" the block in the directions $s, t, u$ and the values are in the big 3-table $R$. We will use the blocks to explain how to fill in the entries of the 3-table and thus to prove that our construction gives indeed *(1) a feasible set of 2-marginals and (2) the entry $R_{1,1,1}$, can be filled in with zero for some 3-table satisfying all 2-marginals if and only if the original set of 2-marginals $A, B, C$ is feasible.* The notation we use in subsequent pictures to depict a way of filling a block is by writing either a single number (e.g. zero), which is used to filled all the block, or by listing a table (e.g. $C_{i,k}$) which indicates the entries of that table are copied down verbatim to be the entries of the block. Figure 6 shows a concrete 3-table which indeed

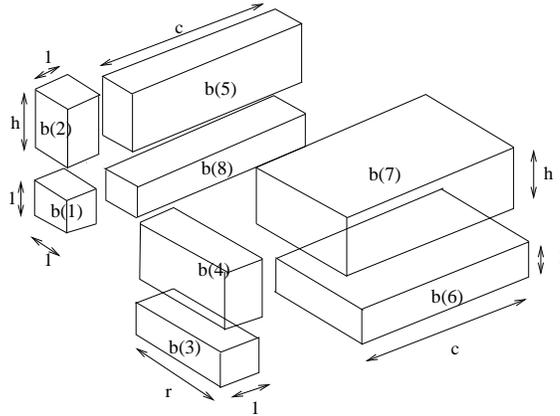

Figure 5: Blocks determined by the proposed 2-marginals.

satisfies the 2-marginals given in the construction. More explicitly we fill the entries as follows:
**Block** $b(1)$, a single entry $R_{1,1,1} = T$. **Block** $b(2)$ has entries $R_{1,1,u}$ for $u = 2, \ldots, h+1$, We fill them by $R_{1,1,u} = 0$. **Block** $b(3)$ has entries $R_{s,1,1}$ for $s = 2, \ldots, r+1$. We fill them by $R_{s,1,1} = 0$. **Block** $b(4)$ has entries $R_{s,1,u}$ for $s = 2..r+1$ and $u = 2, \ldots, h+1$. We fill them by $R_{s,1,u} = C_{s-1,u-1}$. **Block** $b(5)$ has entries $R_{1,t,u}$ for $t = 2..c+1$ and $u = 2, \ldots, h+1$. We fill them by $R_{1,t,u} = B_{t-1,u-1}$. **Block** $b(6)$ has entries $R_{s,t,1}$ for $s = 2..c+1$ and $t = 2, \ldots, r+1$. We fill them by $R_{s,t,1} = A_{s-1,t-1}$. The entries of $b(7)$ and $b(8)$ are all zero. It is simple to verify is that all the axial sums agree with the totals stated in Figure 4 because in the construction we used the 1-marginals of the 2-tables $A, B, C$ as part of the 2-marginals and the data is consistent. This proves the first claim. Now we claim that the entry $R_{1,1,1}$ takes on the value zero for some 3-table $R_{s,t,u}$ of size $(r+1, c+1, h+1)$ if and only if the 2-marginals $A, B, C$ have a feasible solution. Let us assume that there is a 3-table $R_{s,t,u}$ of size $(r+1, c+1, h+1)$ and $R_{1,1,1} = 0$. We divide the argument into two steps illustrated in the left-hand side of Figure 7: Note that if the entry $R_{1,1,1}$ is zero we must have filled $R_{2,1,1} = C_{1,+}, R_{3,1,1} = C_{2,+}, \ldots, R_{(r+1),1,1} = C_{c,+}$ and $R_{1,1,2} = C_{+,1}, \ldots, R_{1,1,2} = C_{+,2}, \ldots R_{1,1,(h+1)} = C_{+,h}$. The reason is 2-marginals $R_{+,1,1}$ and $R_{1,1,+}$ equal the total sum $T$, and $T = \sum C_{i,+} = \sum C_{+,j}$. This completes the filling of blocks 2 and 3. Now, the marginals $R_{s,t,+}$ and $R_{+,t,u}$ attached to $b(4)$ imply that $b(4)$ is simply full of zeros otherwise we surpass the 2-marginals. This completes the first step of the argument.

For the second step we refer to the right of Figure 7. The marginal table $R_{s,+,u}$ and the assignments so far for blocks 2 and 3 imply that only zero values can be put in the entries of blocks

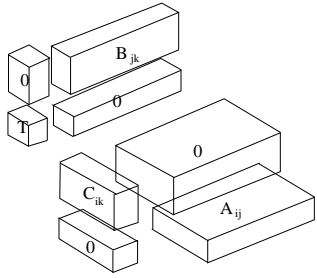

Figure 6: A 3-table with the proposed 2-marginals.

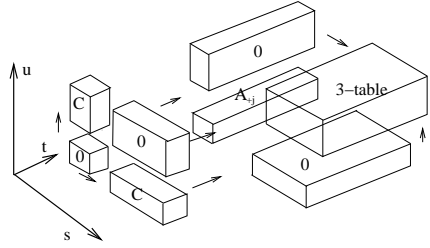

Figure 7: A 3-table with the proposed 2-marginals when $R_{1,1,1} = 0$.

5 and 6, otherwise we surpass $R_{s,+,u}$. Now blocks 7, 8 are left to be decided. The 2-marginals corresponding to $b(8)$ indicate the entries of $b(8)$ are $R_{1,2,1} = A_{+,1} = B_{1,+}, \ldots, R_{1,j,1} = A_{+,j} = B_{j,+}, \ldots R_{1,(c+1),1} = A_{+,c} = B_{c,+}$. Note that these are in fact the 1-marginals that follow from the input 2-marginals $A, B, C$. Finally, block 7 is the only block unfilled. Looking at the zeros that fill block 4, 5, and 6, we see block 7 is indeed a 3-table with 2-marginals $A, B, C$, and this ends the proof of the claim. Now conversely, and essentially following a reverse order, if there is a 3-table with 2-marginals $A, B, C$ we can put a copy of it as block 7. Then by the corresponding 2-marginals we see $b(4), b(5), b(6)$ are filled with zeros. This forces $R_{2,1,1} = C_{1,+}, R_{3,1,1} = C_{2,+}, \ldots, R_{(r+1),1,1} = C_{c,+}$ and $R_{1,1,2} = C_{+,1}, \ldots, R_{1,1,2} = C_{+,2}, \ldots R_{1,1,(h+1)} = C_{+,h}$. This is because the 2 marginals $R_{+,1,1}$ and $R_{1,1,+}$ equal the total sum $T$ and the 2-marginals $R_{s,+,u}$. Finally, the entry $R_{1,1,1}$ is forced to be zero. This completes the proof of the part(A).

Next we prove part (B). We reduce the 3-dimensional matching problem to the problem of deciding whether, given feasible 2-marginals, there is a slim 3-table with a specified entry attaining the Fréchet upper bound. As mentioned in Section 2, the 3-dimensional matching problem is equivalent to the following problem: given a $\{0, 1\}$-valued 3-table $p = (p_{i,j,k})$ of size $(n, n, n)$, is there a 3-table $x = (x_{i,j,k})$ with all 1-marginals $u_{i,+,+}, u_{+,j,+}, u_{+,+,k}$ equal to 1 which is *dominated* by $p$, i.e. satisfies the upper bounds $x_{i,j,k} \leq p_{i,j,k}$ for all $i, j, k$ ? given such data, we expand it to data for upper bounds and 1-marginals for 3-tables of enlarged size $(n+1, n+1, n+1)$ as follows: we maintain the given upper bounds $p_{i,j,k}$ and the 1-marginals $u_{i,+,+}, u_{+,j,+}, u_{+,+,k}$ equal to 1 for $1 \leq i, j, k \leq n$; we introduce the new upper bounds $p_{n+1,n+1,n+1} := 2n$, $p_{i,j,n+1} := p_{i,n+1,k} := p_{n+1,j,k} := 0$ for $1 \leq i, j, k \leq n$, and $p_{i,n+1,n+1} := p_{n+1,n+1,k} := p_{n+1,j,n+1} := 1$ for $1 \leq$ $i, j, k \leq n$. Finally, the three new 1-marginals are introduced by $u_{n+1,+,+} := 2n$, $u_{+,n+1,+} := 2n$, and $u_{+,+,n+1} := 2n$. The extended bounds are shown in Figure 8 on the union of the input $(n, n, n)$-table and seven other blocks. The extended 1-marginals and upper bounds are feasible:

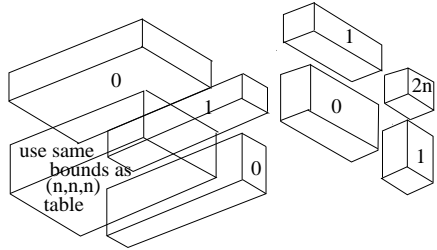

Figure 8: The entry bounds for 3-tables of size $(n+1, n+1, n+1)$.

the $(n+1, n+1, n+1)$-table $x$ defined by setting $x_{i,n+1,n+1} := x_{n+1,n+1,k} := x_{n+1,j,n+1} := 1$ for all $1 \leq i, j, k \leq n$ and zero in all other entries is feasible.

Now consider any feasible extended table $x$: then it is not hard to see that its $(n, n, n)$-subtable $(x_{i,j,k})_{1,1,1}^{n,n,n}$ is feasible for the original data (coming from the input to the 3-dimensional matching) if and only if the entry $x_{n+1,n+1,n+1}$ equals the maximal possible value $2n$.

Now, we "lift" the situation to the problem with 2-marginals and no upper bounds in slim tables as follows: to the extended upper bound and 1-marginal data for $(n+1, n+1, n+1)$-tables, apply the transformation described in Theorem 1.4 (to be proved in the next section). This gives feasible 2-marginals for 3-tables of size $(3, (n+1)^2, 3(n+1))$. By Theorem 1.4, there is a feasible $(n+1, n+1, n+1)$-table $x$ whose entry $x_{n+1,n+1,n+1}$ attains the maximal possible value $2n$ if and only if there is a feasible $(3, (n+1)^2, 3(n+1))$-table $y$ whose entry $y_{1,(n+1)(n+1),\mathrm{dom}\,(n+1)}$ attains the maximal possible value $2n$. This completes the proof of part (B) and the proof of entries $(1,3)$ and hence also entries and $(2,3), (1,4), (2,4)$ in the table of Theorem 1.3.

Finally, we establish entry $(1,2)$ in the statement of Theorem 1.3: we present a polynomial time algorithm for deciding whether there is a 3-table $x$ with specified 2-marginals whose entry $x_{1,1,1}$ is in the range $L \leq x_{1,1,1} \leq U$. The lower (respectively, upper) versions of the entry-security problem is the special cases of this entry-range problem obtained by taking $U$ to be the Fréchet upper bound $U := \min\{v_{1,1,+}, v_{1,+,1}, v_{+,+,1}\}$ (respectively, taking $L := 0$). We use a simple modification of the algorithm for enumeration presented in the proof of Theorem 1.2 in the previous section; using the notation in that proof, we simply need to modify the definition of the first matrix $A_1$, where for $s, t \in S$, its $(s,t)$-th entry is now redefined to be

$$(A_1)_{s,t} := \begin{cases} 1 & \text{if } (t-s)_{i,+} = v_{i,+,k} \text{ for all } i, (t-s)_{+,j} = v_{+,j,k} \text{ for all } j \text{ and } L \leq (t-s)_{1,1} \leq U \\ 0 & \text{otherwise} \end{cases}.$$

The other matrices $A_k$ remain as before. The entry $A_{l,u}^h$ of the product matrix now yields the number of tables with $L \leq x_{1,1,1} \leq U$ and hence is nonzero if and only if such a table exists. □

## 5   Multi-index transportation polytopes and the power of two-marginals

We conclude with the proof of Theorem 1.4 discussed in Subsection 1.4.

**Proof of Theorem 1.4.** The proof is based on an extension of the construction used in the proof of Theorem 1.1. We provide the construction and an abridged form of the argumentation.

Given 1-marginals $(u_{i,+,+})$, $(u_{+,j,+})$, $(u_{+,+,k})$ and entry upper bounds $(p_{i,j,k})$ for 3-tables of size $(r, c, h)$, we define efficiently constructible 2-marginals $(v_{i,j,+})$, $(v_{i,+,k})$, $(v_{+,j,k})$ for 3-tables of size $(3, rc, r + c + h)$ such that nonnegative real arrays $y$ with these marginals are in integer preserving affine bijection with nonnegative real 3-arrays $x$ of size $(r, c, h)$ satisfying the given 1-marginals and upper bounds, thus providing an isomorphism of the corresponding multi-index transportation polytopes and sets of tables of the two systems.

As in the proof of Theorem 1.1, $(3, rc, r + c + h)$-tables will be indexed by triplets with the first index an integer $1 \leq t \leq 3$, the second index an ordered pair $ij$ with $1 \leq i \leq r$ and $1 \leq j \leq c$, and the third index three-letter abbreviation gro $\in \{\text{dom}, \text{row}, \text{col}\}$ along with a numerical index $1 \leq k \leq h$. Let $U$ denote the minimal of the two values $\max\{u_{i,+,+} : 1 \leq i \leq r\}$ and $\max\{u_{+,j,+} : 1 \leq j \leq c\}$. The 2-marginals are provided by the following three matrices.

$$(v_{+,ij,\text{gro}\,k}) = \begin{pmatrix}
 & 11 & 12 & \cdots & 1c & 21 & 22 & \cdots & 2c & \cdots & r1 & r2 & \cdots & rc \\
\text{dom}\,1 & p_{1,1,1} & p_{1,2,1} & \cdots & p_{1,c,1} & p_{2,1,1} & p_{2,2,1} & \cdots & p_{2,c,1} & \cdots & p_{r,1,1} & p_{r,2,1} & \cdots & p_{r,c,1} \\
\text{dom}\,2 & p_{1,1,2} & p_{1,2,2} & \cdots & p_{1,c,2} & p_{2,1,2} & p_{2,2,2} & \cdots & p_{2,c,2} & \cdots & p_{r,1,2} & p_{r,2,2} & \cdots & p_{r,c,2} \\
\cdots & \vdots & \vdots & \vdots & \vdots & \vdots & \vdots & \vdots & \vdots & \vdots & \vdots & \vdots & \vdots & \vdots \\
\text{dom}\,h & p_{1,1,h} & p_{1,2,h} & \cdots & p_{1,c,h} & p_{2,1,h} & p_{2,2,h} & \cdots & p_{2,c,h} & \cdots & p_{r,1,h} & p_{r,2,h} & \cdots & p_{r,c,h} \\
\text{row}\,1 & U & U & \cdots & U & 0 & 0 & \cdots & 0 & \cdots & 0 & 0 & \cdots & 0 \\
\text{row}\,2 & 0 & 0 & \cdots & 0 & U & U & \cdots & U & \cdots & 0 & 0 & \cdots & 0 \\
\cdots & \vdots & \vdots & \vdots & \vdots & \vdots & \vdots & \vdots & \vdots & \vdots & \vdots & \vdots & \vdots & \vdots \\
\text{row}\,r & 0 & 0 & \cdots & 0 & 0 & 0 & \cdots & 0 & \cdots & U & U & \cdots & U \\
\text{col}\,1 & U & 0 & \cdots & 0 & U & 0 & \cdots & 0 & \cdots & U & 0 & \cdots & 0 \\
\text{col}\,2 & 0 & U & \cdots & 0 & 0 & U & \cdots & 0 & \cdots & 0 & U & \cdots & 0 \\
\cdots & \vdots & \vdots & \vdots & \vdots & \vdots & \vdots & \vdots & \vdots & \vdots & \vdots & \vdots & \vdots & \vdots \\
\text{col}\,c & 0 & 0 & \cdots & U & 0 & 0 & \cdots & U & \cdots & 0 & 0 & \cdots & U
\end{pmatrix}$$

$$(v_{t,ij,+}) = \begin{pmatrix}
 & 11 & 12 & \cdots & 1c & 21 & 22 & \cdots & 2c & \cdots & r1 & r2 & \cdots & rc \\
1 & U & U & \cdots & U & U & U & \cdots & U & \cdots & U & U & \cdots & U \\
2 & p_{1,1,+} & p_{1,2,+} & \cdots & p_{1,c,+} & p_{2,1,+} & p_{2,2,+} & \cdots & p_{2,c,+} & \cdots & p_{r,1,+} & p_{r,2,+} & \cdots & p_{r,c,+} \\
1 & U & U & \cdots & U & U & U & \cdots & U & \cdots & U & U & \cdots & U
\end{pmatrix}$$

$$(v_{t,+,\text{gro}\,k}) = \begin{pmatrix}
 & 1 & 2 & 3 \\
\hline
\text{dom}\,1 & u_{+,+,1} & p_{+,+,1} - u_{+,+,1} & 0 \\
\text{dom}\,2 & u_{+,+,2} & p_{+,+,2} - u_{+,+,2} & 0 \\
\ldots & \vdots & \vdots & \vdots \\
\text{dom}\,h & u_{+,+,h} & p_{+,+,h} - u_{+,+,h} & 0 \\
\hline
\text{row}\,1 & c \cdot U - u_{1,+,+} & 0 & u_{1,+,+} \\
\text{row}\,2 & c \cdot U - u_{2,+,+} & 0 & u_{2,+,+} \\
\ldots & \vdots & \vdots & \vdots \\
\text{row}\,r & c \cdot U - u_{r,+,+} & 0 & u_{r,+,+} \\
\hline
\text{col}\,1 & 0 & u_{+,1,+} & r \cdot U - u_{+,1,+} \\
\text{col}\,2 & 0 & u_{+,2,+} & r \cdot U - u_{+,2,+} \\
\ldots & \vdots & \vdots & \vdots \\
\text{col}\,c & 0 & u_{+,c,+} & r \cdot U - u_{+,c,+}
\end{pmatrix}$$

We make use again of a partition of $(3, rc, r + c + h)$-arrays into blocks similar to Figure 1. First consider any nonnegative real $(r, c, h)$-array $x$ satisfying the given 1-marginals and upper bounds. We show it uniquely extends to a nonnegative real $(3, rc, r + c + h)$-array with 2-marginals as above. Given then such an array $x$, embed it in the **black block** $(1, \text{dom})$ of a $(3, rc, r + c + h)$-array $y$ by $y_{1,ij,\text{dom}\,k} := x_{i,j,k}$ for all $i, j, k$. We now show that the block $x$ can be uniquely extended to a whole nonnegative real $(3, rc, r + c + h)$-array $y$ with the above 2-marginals. First, the entries in the **grey blocks** $(1, \text{col})$, $(2, \text{row})$ and $(3, \text{dom})$ in Figure 1 are all zero since so are the 2-marginals $v_{1,+,\text{col}\,k} = v_{2,+,\text{row}\,k} = v_{3,+,\text{dom}\,k} = 0$ for all $k$. Next, consider the entries in the **white block** $(1, \text{row})$: using the fact that all entries in the block $(1, \text{col})$ below it are zero, and examining the 2-marginals $v_{+,ij,\text{row}\,k}$ and $v_{1,ij,+} = U$, we find that $y_{1,ij,\text{row}\,i} = U - \sum_{k=1}^{h} y_{1,ij,\text{dom}\,k} = U - x_{i,j,+} \geq 0$ whereas for $k \neq i$ we have $y_{1,ij,\text{row}\,k} = 0$. This also yields the entries in the **white block** $(3, \text{row})$: we have $y_{3,ij,\text{row}\,i} = U - y_{1,ij,\text{row}\,i} = x_{i,j,+} \geq 0$ whereas for $k \neq i$ we have $y_{3,ij,\text{row}\,k} = 0$. Next, consider the entries in the **white block** $(2, \text{dom})$: using the fact that all entries in the block $(3, \text{dom})$ to its right are zero, and examining the 2-marginals $v_{+,ij,\text{dom}\,k} = p_{i,j,k}$ we find that $y_{2,ij,\text{dom}\,k} = p_{i,j,k} - x_{i,j,k} \geq 0$ for all $i, j, k$. Next consider the entries in the **white block** $(2, \text{col})$: using the fact that all entries in the block $(2, \text{row})$ above it are zero, and examining the 2-marginals $v_{+,ij,\text{col}\,k}$ and $v_{2,ij,+} = p_{i,j,+}$, we find that $y_{2,ij,\text{col}\,j} = p_{i,j,+} - \sum_{k=1}^{h} y_{2,ij,\text{dom}\,k} = x_{i,j,+} \geq 0$ whereas for $k \neq j$ we have $y_{2,ij,\text{col}\,k} = 0$. This also yields the entries in the **white block** $(3, \text{col})$: we have $y_{3,ij,\text{col}\,j} = U - y_{2,ij,\text{col}\,j} = U - x_{i,j,+} \geq 0$ whereas for $k \neq j$ we have $y_{3,ij,\text{col}\,k} = 0$.

Next consider any nonnegative real $(3, rc, r + c + h)$-array $y$ with the above 2-marginals, and let $x$ be its $(r, c, h)$-subarray given by the black block $(1, \text{dom})$ of $y$, defined by $x_{i,j,k} := y_{1,ij,\text{dom}\,k}$ for all $i, j, k$. We show that $x$ is nonnegative and satisfies the given upper bounds and 1-marginals. It is nonnegative since so is $y$. It is dominated by $p$ since, for all $i, j, k$ we have $p_{i,j,k} - x_{i,j,k} = y_{2,ij,\text{dom }k} \geq 0$. Finally, it obeys the 1-marginals $u_{i,+,+}$, $u_{+,j,+}$ and $u_{+,+,k}$ since:

$$x_{+,+,k} = \sum_{i,j} y_{1,ij,\text{dom }k} = v_{1,+,\text{dom }k} = u_{+,+,k}, \quad 1 \leq k \leq h ;$$

$$x_{i,+,+} = \sum_{j} y_{3,ij,\text{row }i} = v_{3,+,\text{row }i} = u_{i,+,+}, \quad 1 \leq i \leq r ;$$

$$x_{+,j,+} = \sum_{i} y_{2,ij,\text{col }j} = v_{2,+,\text{col }j} = u_{+,j,+}, \quad 1 \leq j \leq c .$$

Thus, the set of nonnegative real $(r,c,h)$-arrays $x$ satisfying the given upper bounds and 1-marginals is in integer preserving affine bijection with the set of nonnegative real $(3,rc,r+c+h)$-arrays $y$ with the constructed 2-marginals. In particular, the corresponding multi-index transportation polytopes and sets of tables of the two systems are isomorphic, completing the proof. □

Jesus De Loera

*University of California at Davis, Davis, CA 95616, USA*

*email: jesus@math.ucdavis.edu*

Shmuel Onn

*Technion - Israel Institute of Technology, 32000 Haifa, Israel,*
    *and*
*University of California at Davis, Davis, CA 95616, USA.*
*email: onn@ie.technion.ac.il, onn@math.ucdavis.edu,*
    *http://ie.technion.ac.il/∼onn*